\documentclass[12pt]{article}
\usepackage{amsmath,amssymb,bbm,theorem}
\newcommand{\assign}{:=}
\newcommand{\tmmathbf}[1]{\ensuremath{\boldsymbol{#1}}}
\newcommand{\tmop}[1]{\ensuremath{\operatorname{#1}}}
\newtheorem{theorem}{Theorem}[section]
\newtheorem{corollary}[theorem]{Corollary}
\newtheorem{lemma}[theorem]{Lemma}

\newtheorem{remark}[theorem]{Remark}

\numberwithin{equation}{section}
\newenvironment{proof}{\noindent\textbf{Proof\ }}{\hspace*{\fill}$\Box$\medskip}

\begin{document}

\title{On the volume set of point sets in vector spaces over finite
fields}\author{Le Anh Vinh\\
Mathematics Department\\
Harvard University\\
Cambridge, MA, 02138, US\\
vinh@math.harvard.edu}\maketitle

\begin{abstract}
  We show that if $\mathcal{E}$ is a subset of the $d$-dimensional vector space over a finite field $\mathbbm{F}_q$ ($d \geq 3$) of
  cardinality $|\mathcal{E}| \geq (d-1)q^{d - 1}$, then the set of volumes of
  $d$-dimensional parallelepipeds determined by $\mathcal{E}$ covers
  $\mathbbm{F}_q$. This bound is sharp up to a factor of $(d-1)$ as taking $\mathcal{E}$ to be a $(d -
  1)$-hyperplane through the origin shows.
\end{abstract}

\section{Introduction}

Let $q$ be an odd prime power, and let $\mathbbm{F}_q$ be a finite field of $q$
elements. The distribution of the determinant of matrices with entries in a
restricted subset of $\mathbbm{F}_q$ has been studied recently by various
researchers (see, for example, \cite{as,covert,vinh,vinh2} and the references therein). In particular,
Covert et al. \cite{covert} studied this problem in a more general setting. They examined
the distribution of volumes of $d$-dimensional parallelepipeds determined by a
large subset of $\mathbbm{F}_q^d$. More precisely, for any $\tmmathbf{x}^1,
\ldots, \tmmathbf{x}^d \in \mathbbm{F}_q^d$, define $\tmop{vol}
(\tmmathbf{x}^1, \ldots, \tmmathbf{x}^d)$ as the determinant of the matrix
whose rows are $\tmmathbf{x}^j$s. The focus of \cite{covert} is to study the cardinality
of the volume set
\[ \tmop{vol} (\mathcal{E}) =\{\tmop{vol} (\tmmathbf{x}^1, \ldots,
   \tmmathbf{x}^d) : \tmmathbf{x}^j \in \mathcal{E}\}. \]
A subset $\mathcal{E} \subseteq \mathbbm{F}_q^d$ is called a product-like set
if $|\mathcal{E} \cap \mathcal{H}_n | \ll |\mathcal{E}|^{n / d}$ for any
$n$-dimensional subspace $\mathcal{H}_n \subset \mathbbm{F}_q^d$. Covert et
al. \cite{covert} showed that if $\mathcal{E} \subseteq \mathbbm{F}_q^d$ is a
product-like set of cardinality $|\mathcal{E}| \gg q^{15 / 8}$, then
$\mathbbm{F}_q^{\ast} \subseteq \tmop{vol} (\mathcal{E})$.
When $\mathcal{E} \subseteq \mathbbm{F}_q^3$ is an arbitrary set, they obtained the following
result.

\begin{theorem}
  (\cite[Theorem 2.10]{covert}) Suppose that $\mathcal{E} \subseteq \mathbbm{F}_q^3$ of
  cardinality $|\mathcal{E}| \geqslant C q^2$ for a sufficiently large
  constant $C > 0$. There exists $c > 0$ such that
  \[ | \tmop{vol} (\mathcal{E}) | \geqslant c q. \]
\end{theorem}

In this short note, we show that, under the same condition, $\mathcal{E}$
determines all possible volumes. More precisely, we will prove the following
general result.

\begin{theorem}\label{main-vol}
  When $\mathcal{E} \subseteq \mathbbm{F}_q^d$ and
  $|\mathcal{E}| \geq (d-1) q^{d - 1}$, $\tmop{vol} (\mathcal{E})
  =\mathbbm{F}_q$.
\end{theorem}

\begin{remark}
  The assumption $|\mathcal{E}| \geq (d-1)q^{d - 1}$ is sharp up to a factor of $(d-1)$ as taking
  $\mathcal{E}$ to be a $(d - 1)$-hyperplane through the origin shows.
\end{remark}

Note that the implied constant in the symbol `$\gg$' may depend on
integer parameter $d$. We recall that the notation $U \gg V$ is equivalent to
the assertion that $|U| \geqslant c|V|$ holds for some constant $c > 0$. 

\section{Preparations} \label{sec2-vol}

Recall that
\[ \tmop{vol} (\tmmathbf{x}^1, \ldots, \tmmathbf{x}^d) =\tmmathbf{x}^1 \cdot
   (\tmmathbf{x}^2 \wedge \ldots \wedge \tmmathbf{x}^d), \]
where the dot product is defined by the usual formula
\[ \tmmathbf{u} \cdot \tmmathbf{v}= u_1 v_1 + \ldots + u_d v_d. \]
The generalized cross product, also called the wedge product, is given by
the identity
\[ \tmmathbf{u}^2 \wedge \ldots \wedge \tmmathbf{u}^d = \det
   \left(\begin{array}{c}
     \tmmathbf{i}\\
     \tmmathbf{u}^2\\
     \ldots\\
     \tmmathbf{u}^d
   \end{array}\right), \]
where $\tmmathbf{i}= (\tmmathbf{i}^1, \ldots, \tmmathbf{i}^d)$ indicates the
coordinate directions in $\mathbbm{F}_q^d$.

\subsection{Geometric Incidence Estimates}

One of our main tools is a two-set version of the geometric incidence estimate
developed by D. Hart, A. Iosevich, D. Koh, and M. Rudnev in \cite{hart2} (see also \cite{hart1}
for an earlier version and \cite{covert} for a function version of this estimate). Note
that going from one set formulation in the proof of \cite[Theorem 2.1]{hart2} to a two-set
formulation is just a matter of inserting a different letter into a couple of
places.

\begin{lemma}\label{l1-vol} (\cite[Theorem 2.1]{hart2})
  Let $B (\cdot, \cdot)$ be a non-degenerate bilinear form in
  $\mathbbm{F}_q^d$. For any $\mathcal{E}, \mathcal{F} \subseteq \mathbbm{F}_q^d$, let
  \[ \nu_{t,B} (\mathcal{E}, \mathcal{F}) = \sum_{B(\tmmathbf{x}, \tmmathbf{y})=
     t} \mathcal{E}(\tmmathbf{x})\mathcal{F}(\tmmathbf{y}), \]
\end{lemma}
where, here and throughout the paper, $\mathcal{E}(\tmmathbf{x})$ denotes the
characteristic function of $\mathcal{E}$. We have
\[ \nu_{t,B} (\mathcal{E}, \mathcal{F}) = |\mathcal{E}| |\mathcal{F}|q^{- 1} + R_{t,B}
   (\mathcal{E}, \mathcal{F}), \]
where
\[ |R_{t,B} (\mathcal{E}, \mathcal{F}) |^2 \leqslant |\mathcal{E}|
   |\mathcal{F}| q^{d - 1}, \tmop{if} t \neq 0. \]
As a corollary of Lemma \ref{l1-vol}, D. Hart and A. Iosevich \cite{hart1} derived the following
result.

\begin{corollary}\label{c1-vol}
  (\cite{hart1,hart2}) For any $\mathcal{E}, \mathcal{F} \subseteq \mathbbm{F}_q^d$,
  let
  \[ \mathcal{E} \cdot \mathcal{F}=\{\tmmathbf{u} \cdot \tmmathbf{v}:
     \tmmathbf{u} \in \mathcal{E}, \tmmathbf{v} \in \mathcal{F}\}. \]
  We have $\mathbbm{F}_q^{\ast} \subseteq \mathcal{E} \cdot \mathcal{F}$ when
  $|\mathcal{E}| |\mathcal{F}| \gg q^{d + 1}$.
\end{corollary}

We also need the following corollary. 

\begin{corollary}\label{c2-vol}
  Let $B (\cdot, \cdot)$ be a non-degenerate bilinear form in
  $\mathbbm{F}_q^2$. For any $\mathcal{E} \in \mathbbm{F}_q^2$, let
  \[ B^{\ast} (\mathcal{E}) =\{B (\tmmathbf{x}, \tmmathbf{y}) : \tmmathbf{x},
     \tmmathbf{y} \in \mathcal{E}\} \backslash \{0\}. \]
  We have
  \[ |B^{\ast} (\mathcal{E}, \mathcal{F}) | \geqslant q \left( 1 - \frac{q +
     q^{3 / 2}}{|\mathcal{E}| + q^{3 / 2}} \right) . \]
\end{corollary}

\begin{proof}
  For any $\tmmathbf{x} \in \mathcal{E}$, there exist at most $q$ vectors
  $\tmmathbf{y}$ such that $\tmmathbf{x} \cdot \tmmathbf{y}= 0$. Hence,
  \begin{equation}\label{e1-vol} \sum_{t \in \mathbbm{F}_q^{\ast}} \nu_{t, B} (\mathcal{E}, \mathcal{E})
     \geqslant |\mathcal{E}|^2 - q|\mathcal{E}|. \end{equation}
  Lemma \ref{l1-vol} implies that
  \begin{equation} \label{e2-vol} \nu_{t, B} (\mathcal{E}, \mathcal{E}) \leqslant \frac{|\mathcal{E}|^2}{q}
     + q^{1 / 2} |\mathcal{E}|, \end{equation}
  for any $t \in \mathbbm{F}_q^{\ast}$. The corollary follows immediately from
  (\ref{e1-vol}) and (\ref{e2-vol}).
\end{proof}

\subsection{Cross-product set}

Let $\mathcal{H}$ be any $d$-dimensional vector space over a finite field
$\mathbbm{F}_q$. Let $\{\tmmathbf{v}^1, \ldots, \tmmathbf{v}^d \}$ be an
orthogonal basis of $\mathcal{H}$. For any $d$ vectors $\tmmathbf{u}^1,
\ldots, \tmmathbf{u}^d \in \mathcal{H}$, each vector $\tmmathbf{u}^i$ can be written uniquely as a linear combination of $\{\tmmathbf{v}^1, \ldots, \tmmathbf{v}^d \}$, i.e.
\[ \tmmathbf{u}^i = \sum_{j = 1}^d u_j^i \tmmathbf{v}^j, \,\,\, u^i_j \in \mathbbm{F}_q, 1 \leqslant j \leqslant d.\]
We have
\begin{equation}\label{e3-vol} \tmmathbf{u}^1 \wedge \tmmathbf{u}^2 \wedge \ldots \wedge \tmmathbf{u}^d =
   \det \left( (u^i_j)_{i, j = 1}^d \right) \tmmathbf{v}^1 \wedge \ldots
   \wedge \tmmathbf{v}^d . \end{equation}
For any $\mathcal{E} \subseteq \mathcal{H}$, define
\begin{equation}\label{e4-vol} \mathcal{D}_{\mathcal{E}, d}^{\ast} \assign \left\{ \det \left( (u^i_j)_{i,
   j = 1}^d \right) : \tmmathbf{u}^i = \sum_{j = 1}^d u_j^i \tmmathbf{v}^j \in
   \mathcal{E}, 1 \leqslant i \leqslant d \right\} \backslash \{0\}. \end{equation}
For any $\tmmathbf{x} \in \mathbbm{F}_q^d$ and $\mathcal{E} \subseteq
\mathbbm{F}_q^d$, let
\[ g_{\mathcal{E}} (\tmmathbf{x}) =\#\{(\tmmathbf{u}^1, \ldots,
   \tmmathbf{u}^{d - 1}) \in \mathcal{E}^{d - 1} : \tmmathbf{u}^1 \wedge
   \ldots \wedge \tmmathbf{u}^{d - 1} =\tmmathbf{x}\}. \]
Define the cross-product set of $\mathcal{E}$,
\[ \mathcal{F}_{\mathcal{E}}^{\ast} =\{\tmmathbf{x} \in \mathbbm{F}_q^d :
   g_{\mathcal{E}} (\tmmathbf{x}) \neq 0\} \backslash \{(0, \ldots, 0)\}. \]
For any $\tmmathbf{x} \in \mathbbm{F}_q^d \backslash \{(0, \ldots, 0)\}$, let
$\mathcal{H}^{\tmmathbf{x}} : =\tmmathbf{x}^{\bot} =\{\tmmathbf{y} \in
\mathbbm{F}_q^d : \tmmathbf{x} \cdot \tmmathbf{y}= 0\}$. It is clear that
$\tmmathbf{x} \in \mathcal{F}_{\mathcal{E}} \backslash \{(0, \ldots, 0)\}$ if
and only if there exist $\tmmathbf{u}^1, \ldots, \tmmathbf{u}^{d - 1} \in
\mathcal{H}^{\tmmathbf{x}} \cap \mathcal{E}$ such that
\begin{equation}\label{v-vol} \tmmathbf{u}^1 \wedge \ldots \wedge \tmmathbf{u}^{d - 1} =\tmmathbf{x}. \end{equation}
It follows from (\ref{e3-vol}), (\ref{e4-vol}), and (\ref{v-vol}) that
\[ \mathcal{F}_{\mathcal{E}}^{\ast} \cap \{l\tmmathbf{x}: l \in
   \mathbbm{F}_q^{\ast} \}=\mathcal{D}^{\ast}_{\mathcal{E} \cap
   \mathcal{H}^{\tmmathbf{x}}, d - 1} . \]
Hence, we have proved the following lemma.

\begin{lemma}\label{l2-vol}
  For any $\mathcal{E} \subseteq \mathbbm{F}_q^d$, we have
  \[ |\mathcal{F}_{\mathcal{E}}^{\ast} | = \sum_{H \in G (d - 1, d)} \left|
     \mathcal{D}^{\ast}_{\mathcal{E} \cap H, d - 1} \right|, \]
  where $G (d - 1, d)$ is the set of all $(d - 1)$-dimensional subspaces of
  $\mathbbm{F}_q^d$.
\end{lemma}

\section{Proof of Theorem 1.2} \label{sec3-vol}

The proof proceeds by induction. We first consider the base case, $d = 3$. We
show that if $|\mathcal{E}| > 2 q^2$, then the cross-product set
$\mathcal{F}_{\mathcal{E}}^{\ast}$ is large. From Lemma \ref{l2-vol}, we have
\begin{equation}\label{e5-vol} |\mathcal{F}^{\ast}_{\mathcal{E}} | = \sum_{\mathcal{H} \in G (2, 3)}
   |\mathcal{D}^{\ast}_{\mathcal{E} \cap \mathcal{H}, 2} |. \end{equation}
Since each nonzero vector lies in $(q + 1)$ two-dimensional subspaces of
$\mathbbm{F}_q^2$,
\[ \sum_{\mathcal{H} \in G (2, 3)} |\mathcal{E} \cap \mathcal{H}| = (q + 1)
   |\mathcal{E}|. \]
Let
\[ G^{\mathcal{E}}_{(2, 3)} =\{\mathcal{H} \in G (2, 3) : |\mathcal{E} \cap
   \mathcal{H}| > q\}, \]
we have
\[ \sum_{\mathcal{H} \in G_{(2, 3)}^{\mathcal{E}}} |\mathcal{E} \cap
   \mathcal{H}| > (q + 1) |\mathcal{E}| - q|G (2, 3) | = (q + 1) |\mathcal{E}|
   - q (q^2 + q + 1) > q^3 . \]

Corollary \ref{c2-vol} implies that
\[ |\mathcal{D}^{\ast}_{\mathcal{E} \cap \mathcal{H}, 2} | \geqslant q \left(
   1 - \frac{q + q^{3 / 2}}{|\mathcal{E} \cap \mathcal{H}| + q^{3 / 2}}
   \right), \]
for any $\mathcal{H} \in G (2, 3)$. Since
\[ f (x) = 1 - \frac{q + q^{3 / 2}}{x + q^{3 / 2}} \]
is a concave function on $[q, q^2]$,
\begin{eqnarray}
  \sum_{\mathcal{H} \in G_{(2, 3)}^{\mathcal{E}}}
  |\mathcal{D}^{\ast}_{\mathcal{E} \cap \mathcal{H}, 2} | & \geqslant & q
  \sum_{\mathcal{H} \in G_{(2, 3)}^{\mathcal{E}}} \left( 1 - \frac{q + q^{3 /
  2}}{|\mathcal{E} \cap \mathcal{H}| + q^{3 / 2}} \right) \nonumber\\
  & \geqslant & q \frac{\sum_{\mathcal{H} \in G_{(2, 3)}^{\mathcal{E}}}
  |\mathcal{E} \cap \mathcal{H}|}{q^2} \left( 1 - \frac{q + q^{3 / 2}}{q^2 +
  q^{3 / 2}} \right) \nonumber\\
  & > & q^2 \left( 1 - q^{- 1 / 2} \right) > q^2 / 2. \label{e6-vol}
\end{eqnarray}
It follows from (\ref{e5-vol}) and (\ref{e6-vol}) that $|\mathcal{F}^{\ast}_{\mathcal{E}} | > q^2 /
2$. Hence, $|\mathcal{E}| |\mathcal{F}_{\mathcal{E}}^{\ast} | > q^4$. Corollary \ref{c1-vol}
implies that $\mathbbm{F}_q^{\ast} \subseteq \mathcal{E} \cdot
\mathcal{F}_{\mathcal{E}}^{\ast} \subseteq \tmop{vol} (\mathcal{E})$. By
choosing a matrix of identical rows, we have $0 \in \tmop{vol} (\mathcal{E})$.
The base case $d = 3$ follows.

Suppose that the theorem holds for $d - 1 \geqslant 3$, we show that it also
holds for $d$. Similarly, we show that if $|\mathcal{E}| > (d - 1) q^{d - 1}$
then the cross-product set $\mathcal{F}_{\mathcal{E}}^{\ast}$ is large. Since
each nonzero vector lies in $(q^{d - 1} - 1) / (q - 1)$ $(d - 1)$-dimensional
subspaces of $\mathbbm{F}_q^d$,
\[ \sum_{\mathcal{H} \in G (d - 1, d)} |\mathcal{E} \cap \mathcal{H}| =
   \frac{q^{d - 1} - 1}{q - 1} |\mathcal{E}|. \]
Let
\[ G^{\mathcal{E}}_{(d - 1, d)} =\{\mathcal{H} \in G (d - 1, d) : |\mathcal{E}
   \cap \mathcal{H}| > (d - 2) q^{d - 2} \}, \]
we have
\begin{eqnarray}
  \sum_{\mathcal{H} \in G_{(d - 1, d)}^{\mathcal{E}}} |\mathcal{E} \cap
  \mathcal{H}| & > & \frac{q^{d - 1} - 1}{q - 1} |\mathcal{E}| - (d - 2) q^{d
  - 1} |G (d - 1, d) | \nonumber\\
  & = & \frac{(q^{d - 1} - 1) |\mathcal{E}| - (d - 2) q^{d - 2} (q^d - 1)}{q
  - 1} \nonumber\\
  & > & q^d, \nonumber
\end{eqnarray}
when $q$ is sufficiently large (in fact, $q > d$ is enough). Since
$|\mathcal{E} \cap \mathcal{H}| \leqslant q^{d - 1}$ for each $\mathcal{H} \in
G_{(d - 1, d)}^{\mathcal{E}}$,
\begin{equation} \label{e7-vol} |G_{(d - 1, d)}^{\mathcal{E}} | > q^d / q^{d - 1} = q. \end{equation}
By induction hypothesis, for any $\mathcal{H} \in G_{(d - 1,
d)}^{\mathcal{E}}$,
\begin{equation}\label{e8-vol} |\mathcal{D}^{\ast}_{\mathcal{E} \cap \mathcal{H}, d - 1} | = q - 1. \end{equation}
Putting (\ref{e7-vol}), (\ref{e8-vol}) and Lemma \ref{l2-vol} together, we have
\[ |\mathcal{F}^{\ast}_{\mathcal{E}} | = \sum_{\mathcal{H} \in G (d-1, d)}
   |\mathcal{D}^{\ast}_{\mathcal{E} \cap \mathcal{H}, d-1} | > \sum_{\mathcal{H}
   \in G_{(d-1, d)}^{\mathcal{E}}} |\mathcal{D}^{\ast}_{\mathcal{E} \cap
   \mathcal{H}, d-1} | > q (q - 1) > q^2 / 2. \]
Hence, $|\mathcal{E}| |\mathcal{F}_{\mathcal{E}}^{\ast} | > q^{d + 1}$. Corollary \ref{c1-vol} implies that $\mathbbm{F}_q^{\ast} \subseteq \mathcal{E} \cdot
\mathcal{F}_{\mathcal{E}}^{\ast} \subseteq \tmop{vol} (\mathcal{E})$. By
choosing a matrix of identical rows, we have $0 \in \tmop{vol} (\mathcal{E})$.
This completes the proof of the theorem.


\begin{thebibliography}{99}

\bibitem{as} O. Ahmadi and I. E. Shparlinski, Distribution of matrices with restricted entries over finite fields, \textit{Inda. Mathem.} \textbf{18}(3) (2007), 327--337.

\bibitem{covert} D. Covert, D. Hart, A. Iosevich, D. Koh, and M. Rudnev, Generalized incidence theorems, homogeneous forms and sum-product estimates in finite fields, \textit{European Journal of Combinatorics}, to appear. 

\bibitem{hart1} D. Hart and A. Iosevich, Sums and products in finite fields: an integral geometric viewpoint, \textit{Contemporary Mathematics}, Radon transforms, geometry, and wavelets, \textbf{464} (2008).


\bibitem{hart2} D. Hart, A. Iosevich, D. Koh, and M. Rudnev, Averages over hyperplanes, sum-product theory in vector spaces over finite fields and the Erd\"os-Falconer distance conjecture, preprint (2007).

\bibitem{vinh} L. A. Vinh, On the distribution of determinants of matrices with restricted entries over finite fields, \textit{Journal of Combinatorics and Number Theory}, to appear.

\bibitem{vinh2} L. A. Vinh, Singular matrices with restricted rows in vector spaces over finite fields, submitted.

\end{thebibliography}
\end{document}